\documentclass{article}

\usepackage{amssymb}
\usepackage{amsmath}

\newtheorem{theorem}{Theorem}
\newtheorem{lemma}{Lemma}[section]
\newtheorem{propo}{Proposition}[section]

\newcommand{\ep}{\varepsilon}
\newcommand{\la}{\lambda}

\newcommand{\R}{\mathbb{R}}

\def\limn{\lim_{n\to\infty}}
\def\lime{\lim_{\ep\to 0}}
\def\Ker{\mbox{Ker}}
\def\Im{\mbox{Im}}
\def\osb{\overline{\sigma}}
\def\usb{\underline{\sigma}}
\def\sn{\sigma_n}
\def\e{\ep}

\def\proof{\smallskip\noindent{\it Proof.} }
\newcommand{\qed} {\hspace {0.1in} \rule {1.5mm} {3.5mm}}
\title{Betti numbers are testable \footnote
{Research sponsored by OTKA Grant K 69062 and 49841 }}

\author{G\'abor Elek}

\begin{document}
\maketitle

\begin{abstract}
We prove that the Betti numbers of simplicial complexes of bounded 
vertex degrees are testable in constant time.
\end{abstract}

\section{Introduction}
Property testing in bounded degree graphs was introduced in the paper
of Goldreich and Ron \cite{GR}. 
In this paper we study property testing for bounded degree simplicial
complexes in higher dimensions.  
 Let $d\geq 2$ be
a natural number and consider finite simplicial complexes
where each vertex (zero dimensional simplex) is contained in at most
$d$ edges ($1$-dimensional simplex). Of course, such a complex can be at
 most $d$-dimensional. What does it mean to test the
$p$-th Betti number of such a simplicial complex ?
First fix a positive real number $\ep>0$. A tester takes a simplicial 
complex $K$ as an input and pick $C(\ep)$ random vertices. Then it
looks at the $C(\ep)$-neighborhoods of the chosen vertices. Based
on this information the tester gives us
 a guess $\hat{b}^p(K)$ for the $p$-th Betti number
$b^p(K)$ of the simplicial complex such a way that :
$$\mbox{Prob}\,(\frac{|\hat{b}^p(K)-b^p(K)|}{|V(K)|}>\ep)<\ep\,,$$
where $V(K)$ is the set of vertices in $K$.
In other words, we can estimate the $p$-th Betti number very effectively
with high probability knowing only a small (random) part of the
simplicial complex.
The goal of this paper is to show the existence of such a tester
for any $\ep>0$. That is to prove the following theorem:
\begin{theorem} \label{tetel1}
Betti-numbers are testable for bounded degree simplicial complexes.
\end{theorem}
For graphs the $0$-th Betti number is just the number of components
and the first Betti number can be computed via the $0$-th Betti number
and the Euler-characteristic, hence it is not hard to see that
such tester exists.
For connected surfaces one can also calculate the first Betti number using
just the number of vertices, edges and triangles. However in higher dimensions
there is no such formula even for triangulated manifolds. 
Note that this paper was not solely
motivated by the paper of Goldreich and Ron, but also by the solution of
the Kazhdan-Gromov Conjecture by Wolfgang L\"uck \cite{Lueck1}
. The workhorse lemma
of our paper is basically extracted from his paper using a slightly
different language. It is very important to note that
our proof works only for Betti numbers of real coefficients and
we do not claim anything for the Betti numbers of
$\mbox{mod-$p$}$\,coefficients. 
\section{The convergence of simplicial complexes}
Let $\Sigma^d$ be the set of finite simplicial complexes $K$
of vertex degree bound $d$ that is any $0$-dimensional simplex
is contained in at most $d$ $1$-dimensional simplices. We denote by $K_i$
the set of $i$-simplices in $K$ and by $G_K$ the $1$-skeleton of $K$, that
is $V(G_K)=K_0$, $E(G_K)=K_1$. A rooted $r$-ball of degree bound $d$ is
a simplex $L\in\Sigma^d$ with a distinguished vertex $x$ such that
for any $y\in V(G_L)$, $d(x,y)\leq r$, where $d(x,y)$ is the
shortest path distance of $x$ and $y$ in the graph $G_L$.
We denote by $Z^{r,d}$ the rooted isomorphism classes of rooted $r$-balls.
If $K\in\Sigma^d$ and $p\in V(K)$ then let $G_r(p)$ be the rooted $r$-ball in
the $1$-skeleton $G_K$ and $B_r(p)$ is the set of simplices $\sigma$ such that
all vertices of $\sigma$ are in $G_r(p)$. Then $B_r(p)$ is a rooted $r$-ball of
vertex degree bound $d$.

\noindent
For $\alpha\in Z^{r,d}$ we denote by $T(K,\alpha)$ the set of vertices $p$
such that $B_r(p)\cong\alpha$. We set
$$p_K(\alpha):=\frac{|T(K,\alpha)|}{|K_0|}\,.$$
We say that $\{K^n\}^\infty_{n=1}\subset \Sigma^d$
is convergent (see \cite{BS} for the graph case) if
\begin{itemize}
\item
$|K^n_0|\to \infty$
\item $\lim_{n\to\infty} p_{K_n}(\alpha)$ exists for
any $r\geq 1$ and $\alpha\in Z^{r,d}$.
\end{itemize}
It is easy to see that any sequence $\{L^n\}^\infty_{n=1}\subset \Sigma^d$
such that $|L^n_0|\to \infty$ contains a convergent subsequence.

\noindent 
Let $\alpha_1, \alpha_2, \alpha_3,\dots$ be an enumeration of 
all the $r$-balls, $r\geq 1$. Then we have the pseudo-metric
$$d_s(K,L):=\sum^\infty_{i=1}\frac{1}{2^i} |p_K(\alpha_i)-p_L(\alpha_i)|\,.$$
Clearly, if $\{K^n\}^\infty_{n=1}$ be an increasing sequence of
simplicial complexes they are Cauchy if and only if they are convergent.
By an oriented $d$-complex $Q$ we mean an element of 
$\Sigma^d$ equipped with a fixed orientation for each of its simplex. Note
that we do not assume that the orientations are compatible in any sense.
We denote by $\hat{\Sigma}^d$ the set of all finite oriented $d$-complexes.
We also define oriented $r$-balls, the set $\hat{Z}^{r,d}$ of all
oriented $r$-ball isomorphism classes and the probabilities
$p_Q(\beta)$ accordingly. Naturally we can define the convergence of
oriented $d$-complexes as well.
\begin{propo}
Let $\{K^n\}^\infty_{n=1}\subset \Sigma^d$ be a convergent
sequence of $d$-complexes. Then one has an oriented copy $Q^n$
for each $K^n$ such that
$\{Q^n\}^\infty_{n=1}\subset \hat{\Sigma}^d$ is convergent as well.
\end{propo}
\proof
Consider i.i.d  random variables
$\theta(x)$ distributed uniformly on $[0,1]$.
Let $(a_0,a_1,\dots,a_i)\in K^n_i$ be a positive orientation if 
$$\theta(a_0)<\theta(a_1)<\dots <\theta(a_i)\,.$$
Then the resulting sequence $\{Q^n\}^\infty_{n=1}\subset \hat{\Sigma}^d$
is convergent with probability $1$. For details of the simple argument see
(\cite{ElekRSA}, Proposition 2.2)\qed
\vskip0.2in
Finally, we need a technical definition that we use in the subsequent
sections. 
Let us consider the set $Q_i$ of $i$-simplices in an oriented $d$-complex $Q$.
If $\sigma,\tau\in Q_i$ we say that $\sigma$ and $\tau$ are adjacent
if they have at least one joint vertex.
This way we can define the shortest path distance $d_i(\sigma,\tau)$
(where $d_i(\sigma,\tau)=\infty$ is possible). The ball $B^i_r(\sigma)$
is the set of $i$-simplices $\tau$ such that $d_i(\sigma,\tau)\leq r$.
As above, we can define the classes $\hat{Z}_i^{r,d}$ and the
sampling probabilities $p_{Q_i}(\beta)$. Note that if $\alpha,\beta\in
\hat{Z}_i^{r,d}$ we say that $\alpha$ is isomorphic to $\beta$ if
they are isomorphic as simplicial complexes not only as metric spaces.

That is $p_{Q_i}(\beta)$ defined the following way.
Let $\beta\in\hat{Z}_i^{r,d}$. Denote by $T(Q_i,\beta)$ the number of
$i$-simplices $\tau$ such that the simplicial complex $B^i_r(\sigma)$
is isomorphic to $\beta$, where the isomorphism preserves the root-simplex.
Then
$$p_{Q_i}(\beta):=\frac{|T(Q_i,\beta)|}{|Q_i|}\,.$$

Clearly, if $\{Q^n\}^\infty_{n=1}\subset \hat{\Sigma}^d$
is a convergent sequence then for any $\beta\in \hat{Z}_i^{r,d}$,
$\limn p_{Q^n_i}(\beta)$ exists.
\section{Betti numbers and combinatorial Laplacians}
Let $Q\in\hat{\Sigma}^d$ be an oriented
simplicial complex.
Let $C^i(Q)$ denote the euclidean space of real functions on
the $1$-simplices of $Q$. Let us consider the cochain-complex
$$C^0(Q)\stackrel {d_0}{\to} C^1(Q)\stackrel {d_1}{\to}\dots$$
Recall that if $f\in C^q(Q)$ then 
$$df(a_0, a_1,\dots, a_{q+1})=f(a_1, a_2,\dots, a_q)- f(a_0,a_2,\dots,
a_{q+1})+\dots$$$$+(-1)^{q+1}f(a_0,a_1,\dots,a_q)\,.$$
Then $b^i(Q)=\dim \Ker\,d_i-\dim \Im\, d_{i-1}$
are the Betti numbers of $Q$. Note that they do not depend on
the choice of the orientation of $Q$ only the underlying simplicial complex.

\noindent
The combinatorial Laplacians (see e.g. \cite{ElekLap})
$\Delta^i_Q:C^i(Q)\to C^i(Q)$ are defined as
$$\Delta^i_Q:=d_{i-1}d^*_{i-1}+d_i^*d_i\,.$$
The operators $\Delta^i_Q$ are positive and self-adjoint. Also,
$$\dim\Ker\Delta^i_Q=b^i(Q)\,.$$
Let us remark that by Lemma 2.5 of \cite{ElekLap}
we have the following information on the combinatorial Laplacians:
\begin{itemize}
\item 
$\Delta^i_Q(\sigma,\tau)\neq 0$ only if $\sigma=\tau$ or $\sigma$ and $\tau$
are adjacent.
\item 
$\Delta^i_Q(\sigma,\tau)$ is always an integer.
\item
$|\Delta^i_Q(\sigma,\tau)|\leq d+1$.
\end{itemize}
\section{Weak convergence of probability measures}

First recall the notion of weak convergence of probability measures.
Let $\{\mu_n\}^\infty_{n=1}$ be probability measures
on the interval $[0,K]$.
Then $\{\mu_n\}^\infty_{n=1}$  weakly converges to $\mu$ if for any
continuous function $f\in C[0,K]$
$$\int^K_0 f d\mu_n\to \int^K_0 f d\mu\,.$$
For an example, let $\mu_n(\frac{1}{n})=1$, then
the measures $\mu_n$ converge to the measure concentrated at the zero.
Note that in this case $\limn \mu_n(0)\neq\mu(0)\,.$

\noindent
Now let $$c(\mu_n)=
\lime \int_\ep^K \log \lambda d\mu_n\,.$$

The following theorem can be extracted from \cite{Lueck1},
nevertheless we provide a proof using only the language of real
analysis, avoiding any reference to operators.
\begin{theorem}
Suppose that $\{\mu_n\}^\infty_{n=1}$  weakly converges to $\mu$
and for any $n\geq 1$, $c(\mu_n)\geq 0\,.$ Then
$\limn \mu_n(0)=\mu(0)\,.$ \end{theorem}
\proof
First we need some notations.
For a monotone function $f$,
$$f^+(\lambda)=\inf_{\ep\to 0} f(\lambda+\ep)\,$$
For the measures $\mu_n$ let $\sigma_n$ be their distribution
function that is
$$\sigma_n(\lambda)=\mu_n([0,\lambda])\,.$$
Also, let $t(\lambda)=\mu([0,\lambda])\,.$
Note that $\sigma^+_n=\sigma_n\,.$
Let
$$\osb(\lambda):=\limsup_{n\to\infty} \sn(\lambda)$$
and
$$\usb(\lambda):=\liminf_{n\to\infty} \sn(\lambda)\,.$$
The following lemma trivially follows from
the definitions.
\begin{lemma}\label{l73}
Let $f$ be a continuous function
such that
$$\chi_{[0,\lambda]}(x)\leq f(x)\leq 
\chi_{[0,\lambda+\frac{1}{k}]}(x)+ \frac{1}{k}$$ for any $0\leq x \leq K$,
then
$$\sigma_n(\lambda)\leq \int_0^K f(\lambda) d\mu_n \leq
\sigma_n(\lambda+\frac{1}{k})+\frac{1}{k}\,.$$
$$t(\lambda)\leq \int_0^K f(\lambda) d\mu \leq
t(\lambda+\frac{1}{k})+\frac{1}{k}\,.$$ \end{lemma}
\begin{propo} \label{p74}
$\osb(\lambda)\leq t(\lambda)=\usb^+(\lambda)=\osb^+(\lambda)\,.$ \end{propo}
\proof
By this lemma,
$$\osb(\lambda)\leq t(\lambda)\leq\osb(\lambda+\frac{1}{k})+\frac{1}{k}\,.$$
Hence $\osb(\lambda)\leq t(\lambda)=\usb^+(\lambda)\leq \osb^+(\lambda)\,.$
Since $t(\lambda)$ is monotone,
we have that
$\osb(\lambda+\ep)\leq t(\lambda+\epsilon)$ and
$\osb^+(\lambda)\leq t^+(\lambda)=t(\lambda)\,.$
Thus our proposition follows.\qed

\vskip0.2in
\noindent
The following elementary analysis lemma is proved in \cite{Lueck2}.
\begin{lemma} \label{stieltjes}
Let $f$ be a continuously differentiable function on the positive reals
and $\mu$ be a probability measure on the $[0,K]$ interval. Suppose
that $F$ is the distribution function of $\mu$, that is $\mu[0,\lambda]=
F(\lambda)$. Then for any $0<\ep\leq K$:
$$\int_\e^K f(\lambda) d\mu= -\int^K_\e f'(\lambda)F(\lambda) d\lambda
+ f(K)F(K)-f(\e)F(\e)\,.$$
\end{lemma}
Assume that $K\geq 1$, then by the previous lemma
$$c(\mu_n)=\log K-\log(\e)\sigma_n(\e)-\int_\e^K\sn(\la)\frac{1}{\la} d\la\,.$$
Since
$$\int^K_\e\frac{\sn(0)}{\la}=(\log K-\log(\e))\sn(0)\,,$$
we have that
$$c(\mu_n)=\log K (1-\sn(0))-\int^K_0\frac{\sn(\la)-\sn(0)}{\la}\,d\la\,.$$
That is 
\begin{equation}
\label{e71}
\int^K_0\frac{\sn(\la)-\sn(0)}{\la}\leq \log K.
\end{equation}
Observe that
$$\int^K_\e\frac{\usb(\la)-\osb(0)}{\la} d\la\leq
\int^K_\e \frac{\liminf_{n\to\infty}(\sn(\la)-\sn(0))}{\la} d\la\,.$$
By Fatou's Lemma,
\begin{equation} \label{e72}
\int^K_\e \frac{\liminf_{n\to\infty}(\sn(\la)-\sn(0))}{\la} d\la\leq
\liminf_{n\to\infty} 
\int^K_\e \frac{\sn(\la)-\sn(0)}{\la} d\la\,.
\end{equation}
Since the right hand side of (\ref{e72}) is less than $\log K$, we
obtain the following inequality:
$$\int_0^K\frac{\usb(\la)-t(0)}{\la}\,d\la\leq \log K\,.$$
Therefore,
$\lim_{\la\to 0} \usb(\lambda)=\osb(0)\,.$
That is by Proposition \ref{p74}
\begin{equation} \label{e73}
\osb(0)=t(0).
\end{equation}
Since one can apply (\ref{e73}) for any subsequence of $\{\sn\}^\infty_{n=1}$
we obtain that
$$\limn \sn(0)=t(0)\,.
\quad\qed$$

\section{Spectral convergence}
The goal of this section is to prove the main technical
proposition of our paper. Note this is based again on the ideas in
\cite{Lueck1}.
Let $P:\R^n\to\R^n$ be a positive, self-adjoint operator and
$\mu_P$ be its normalized spectral measure that is
$$\mu_P(\lambda):=\frac{\mbox{the multiplicity of $\la$ as an
eigenvalue}}{n}\,.$$
Note that if $\|P\|\leq K$ then $\mu_P$ is concentrated on
the interval $[0,K]$.
\begin{propo} \label{spectral}
Let $\{Q^n\}^\infty_{n=1}\subset \hat{\Sigma}^d$ be a convergent
sequence of oriented simplicial complexes. Then there exists $K>0$
such that
\begin{enumerate}
\item
For any $i\geq 1$ and $n\geq 1$, $\|\Delta_{Q^n}^i\|\leq K\,.$
\item
The normalized spectral measures
of $\{\Delta_{Q^n}^i\}^\infty_{n=1}$, $\{\mu^i_n\}^\infty_{n=1}$ 
weakly converge.
\item
$c(\mu^i_n)\geq 0$ for any $i,n\geq 1$.
\end{enumerate}
\end{propo}
\proof
To show (a) it is enough to prove the following lemma.
\begin{lemma}
Let $L,M>0$ be positive integers. Then if $A$ is a $n\times n$-matrix
of real coefficients (that is a linear operator on $\R^n$) such that
\begin{itemize}
\item each row and column of $A$ contains at most $L$ non-zero elements
\item for each entry $A_{i,j}$, $|A_{i,j}|\leq M$ 
\end{itemize}
then $\|A\|\leq 2LM.$
\end{lemma}
\proof
For unit-vectors $f,g\in \R^n$
$$|\langle A(f),g\rangle|=
|\sum_{1\leq i,j\leq n} A_{i,j}f(i)g(j)|\leq
M\sum_{1\leq i,j\leq n}| f(i)g(j)|\,.$$
That is
$$|\langle A(f),g\rangle|\leq \sum_{1\leq i,j\leq n} (f(i)^2 +g(j)^2)\,.$$
Note that the number of occurences of each $f(i)^2$ or $g(j)^2$ is at most $L$
hence $|\langle A(f),g\rangle|\leq 2LM$.\qed

\noindent
Now let us turn to part (b).
The convergence of $\{\mu^i_n\}^\infty_{n=1}$ means that
$$\limn\int^K_0 P(t)\,d\mu^i_n(t)$$
exists for any real polynomial $P$.
That is one needs to prove that
$$\limn\frac{\sum^{|Q^n_i|}_{j=1} (\lambda^{i,n}_j)^r}{|Q^n_i|}$$
exists where
$\{\lambda^{i,n}_j\}$ denotes the spectrum of the $i$-th Laplacian of $Q^n$. 
Hence it is enough to
prove that the limit of normalized traces
\begin{equation}\label{vege}
\limn\frac{\sum_{\sigma\in Q^n_i}(\Delta^i_{Q^n}
  (\sigma,\sigma))^r}{|Q^n_i|}\end{equation}
exists.
The value of $\Delta^i_{Q^n}
  (\sigma,\sigma)$ depends only on the $r$-neighboorhood of $\sigma$,
  therefore the convergence of the complexes $\{Q^n\}^\infty_{n=1}$
immediately implies the existence of the limit in (\ref{vege}).

\noindent
Part (c) follows from the simple fact: If $Q$ is a symmetric integer
matrix then the product of its non-negative eigenvalues is an integer
as well. Indeed, let $\lambda_1,\lambda_2,\dots,\lambda_q$ the list
of the non-zero eigenvalues of $Q$ with multiplicities.
Let $p(t)=det(tI-Q)$ be the characteristic polynomial of $Q$.
Then $p(t)=t^sq(t)$, where $q(0)\neq 0$. Obviously, $q$ is an integer
polynomial, and $|q(0)|=|\prod^q_{i=1} \lambda_q|$. \qed
\section{The proof of Theorem \ref{tetel1}}
We need to prove the following lemma.
\begin{lemma} \label{l61}
Let $\{K^n\}^\infty_{n=1}\subset\Sigma^d$ be simplicial complexes, then
$\limn \frac{b^i(K_n)}{|V(K^n)|}$ exists for any $i\geq 1$.
\end{lemma}
Indeed,  by Lemma \ref{l61} we can immediately see that for any $\e>0$
there exists $\delta>0$ such that if $d_s(K,L)<\delta$ for some $K,L\in
\Sigma^d$, then
\begin{equation}
\label{e61}
\left|\frac{b^i(K)}{|V(K)|}-\frac{b^i(L)}{|V(L)|}\right|\leq\e\,,
\end{equation}
for any $i\geq 1$.
By the definition of the metric $d_s$ there exists some $r\geq 1$
and $\rho>0$ such that if
$$|p_M(\alpha)-p_N\alpha)|\leq \rho$$
for any $r'\leq r$ and $\alpha\in Z^{r',d}$ then 
$d_s(M,N)<\delta$ for any $M,N\in\Sigma^d$.

\noindent
Now let $L^1,L^2,\dots,L^m$ be a finite set of simplicial
complexes such that for any finite simplicial complex $K$ there
exists $1\leq j\leq m$ such that
\begin{equation}
\label{e62}
|p_K(\alpha)-p_{L^j}(\alpha)|\leq\frac{\rho}{10}
\end{equation}
for any  $r'\leq r$ and $\alpha\in Z^{r',d}$, where $r,\rho$ are the
constants above.  The existence of such finite system is clear from
compactness.

\noindent
By the classical Chernoff's inequality there exists $N_\epsilon>0$ such that
the following holds:

\noindent
Let $M\in\Sigma^d$ be an arbitrary simplicial complex. Pick $N_\e$ random
vertices of $M$ and for any  $r'\leq r$ and $\alpha\in Z^{r',d}$
let $Q(M,\alpha)$ be the number of picked vertices $x$ such that
$B_{r'}(x)\cong\alpha$. Then
\begin{equation}
\label{e63}
\mbox{Prob}\left\{ \left|\frac{Q(M,\alpha)}{N_\e}-p_M(\alpha)
\right|>\frac{\rho}{10}\,
\mbox{for at least one $\alpha$}\right\}<\epsilon
\end{equation}
Thus we have the following testing algorithm.
Take the simplicial complex $M$ as an input. Pick $N_\e$ random vertices
and calculate $Q(M,\alpha)$ for all $r'\leq r$, $\alpha\in Z^{r,d}$,
where $r$ is the constant above. Check the list $L^1, L^2,\dots, L^m$.
By (\ref{e63}), with probability more than $(1-\e)$ we find an $L^j$ such that
$|Q(M,\alpha)-P_{L_j}(\alpha)|<\frac{\rho}{5}$
for any $\alpha$. Let $b^i(L_j)$ be our guess. Then by (\ref{e61}) with
probability more than $1-\e$
$$\left|\frac{b^i(L^j)}{|V(L^j)|}- \frac{b^i(M)}{|V(M)|}\right|<\epsilon\,.$$

So, let us prove Lemma \ref{l61}.
By Proposition \ref{spectral}
$$\limn \frac{\dim\Ker \Delta^i_{K^n}}{|K^n_i|}$$
exists. Also, by the definition  of the convergence of simplicial complexes
$\limn \frac{|K^n_i|}{|V(K^n_i)|}$ exists, hence the lemma, and thus our 
Theorem
follows. \qed

\end{document}